\documentclass[12pt]{amsart}
\newtheorem{thm}{Theorem}

\newcommand{\spn}{\operatorname{span}}

\begin{document}
\title[Grushin and Martinet distributions]{Integrability conditions for the
Grushin and Martinet distributions}
\author[Ovidiu Calin]{Ovidiu Calin}
\address{\hskip-\parindent
Department of Mathematics, Eastern Michigan University,\newline
Ann Arbor, MI 48197, USA}
\email{ocalin@emich.edu}
\author[Der-Chen Chang]{Der-Chen Chang}
\address{\hskip-\parindent
Department of Mathematics and Statistics,\newline
Georgetown University, Washington, DC 20057, USA\newline
Department of Mathematics, Fu-Jen Catholic University,\newline
Taipei, Taiwan 24205, ROC}
\email{chang@georgetown.edu}
\author[Michael Eastwood]{Michael Eastwood}
\address{\hskip-\parindent
Mathematical Sciences Institute, Australian National University,\newline
ACT 0200, Australia}
\email{meastwoo@member.ams.org}
\footnotetext{The second author is partially supported by an
NSF grant DMS-1203845 and a Hong Kong RGC competitive earmarked
research grant $\#$601410. 
The third author is supported by the Australian Research Council.}
\begin{abstract}
We realise the first and second Grushin distributions as symmetry reductions of
the 3-dimensional Heisenberg distribution and 4-dimensional Engel distribution
respectively. Similarly, we realise the Martinet distribution as an
alternative symmetry reduction of the Engel distribution. These reductions
allow us to derive the integrability conditions for the Grushin and Martinet
distributions and build certain complexes of differential operators. These
complexes are well-behaved despite the distributions they resolve being
non-regular.
\end{abstract}
\renewcommand{\subjclassname}{\textup{2010} Mathematics Subject Classification}
\maketitle

\section{Introduction}
For each $k\geq 0$, the pair of linear differential operators on 
${\mathbb{R}}^2$ 
$$X\equiv \partial/\partial x\qquad Y\equiv x^k\partial/\partial y$$
generate what is known as a {\em Grushin distribution\/}~\cite{G}. For
$k\geq 1$, it is not really a distribution in the classical sense because the
span of $X$ and $Y$ drops rank along the $y$-axis, $\{x=0\}$. Nevertheless, the
fields $X$ and $Y$ are {\em bracket generating\/} in the sense that taking 
sufficiently many Lie brackets amongst them generates all vector fields. For 
example, if $k=1$, then 
$$X\qquad Y\qquad Z\equiv[X,Y]=\partial/\partial y$$
span all vector fields at which point we notice that
$$[X,Z]=0\qquad[Y,Z]=0.$$
Similarly, if $k=2$, then 
$$X\qquad Y\qquad Z\equiv[X,Y]=2x\partial/\partial y
\qquad W\equiv[X,Z]=2\partial/\partial y$$
span all vector fields at which point all other commutators between these
fields vanish. In this article, we shall only be concerned with the cases 
$k=1$ and $k=2$ but, in fact, the case $k=0$ is very familiar for then the 
integrability equation for the system
$$\left.\begin{array}{crcl}Xf&=&a\\ 
Yf&=&b\end{array}\right\}\quad\mbox{is}\enskip Xb=Ya.$$
More precisely, if we denote by ${\mathcal{E}}$ the germs of smooth functions 
on~${\mathbb{R}}^2$, then the complex of differential operators
$$\begin{picture}(160,40)
\put(0,20){\makebox(0,0){$0$}}
\put(10,20){\vector(1,0){20}}
\put(40,20){\makebox(0,0){${\mathcal{E}}$}}
\put(50,22){\vector(2,1){20}}
\put(50,18){\vector(2,-1){20}}
\put(80,34){\makebox(0,0){${\mathcal{E}}$}}
\put(80,20){\makebox(0,0){$\oplus$}}
\put(80,6){\makebox(0,0){${\mathcal{E}}$}}
\put(90,32){\vector(2,-1){20}}
\put(90,8){\vector(2,1){20}}
\put(120,20){\makebox(0,0){${\mathcal{E}}$}}
\put(130,20){\vector(1,0){20}}
\put(160,20){\makebox(0,0){$0$}}
\end{picture}$$
is locally exact except at the leftmost ${\mathcal{E}}$ where the cohomology 
is~${\mathbb{R}}$, the real-valued locally constant functions. This is the 
familiar de~Rham complex with local exactness a consequence of the Poincar\'e 
Lemma. 

The aim of this article is to present similar integrability conditions and
consequent differential complexes for the Grushin distribution in cases
$k=1,2$ and also for the Martinet distribution~\cite{M}, which is a pair of 
differential operators on ${\mathbb{R}}^3$ as follows:
\begin{equation}\label{martinet}
X=\partial/\partial x\qquad Y=\partial/\partial z + x^2\partial/\partial y.
\end{equation}
This is not regular in the classical sense because the span of the derived
vector fields 
$$X\qquad Y\qquad Z\equiv[X,Y]=2x\partial/\partial y$$
drops rank along the $(y,z)$-plane, $\{x=0\}$. But, again, the
two fields $X$ and $Y$ are bracket generating since
$$X\qquad Y\qquad Z\qquad W\equiv[X,Z]=2\partial/\partial y$$
span the full tangent space.

The complexes that we shall construct are not true resolutions. In addition to
having cohomology equal to ${\mathbb{R}}$ at the start, we shall allow them to
have finite-dimensional cohomology in other degrees. Bearing in mind the lack
of regularity, i.e.~that $\spn\{X,Y\}$ and other vector spaces generated by Lie
bracket are allowed to jump in dimension from point to point, this is a small
price to pay.

We shall devote separate sections of this article to the three cases under
consideration and consign to two appendices brief reviews of the Heisenberg and
Engel distributions from which they will be derived.

This work was initiated while the authors were participating in the
International Workshop on Several Complex Variables and Complex Geometry,
which was held during July 9-13, 2012 at the Institute of Mathematics, Academia
Sinica, Taipei. The authors would like to thank the local organizers,
especially Professor Jih-Hsin Cheng for the invitation and warm hospitality
extended to them during their visit to Taiwan.

\section{The first Grushin distribution}\label{firstGrushin}
Recall that in this case we are concerned with the three vector fields
\begin{equation}\label{firstgrushin}
X=\partial/\partial x\qquad Y=x\partial/\partial y\qquad 
Z\equiv[X,Y]=\partial/\partial y\end{equation}
on ${\mathbb{R}}^2$ with co\"ordinates~$(x,y)$. We shall derive integrability
conditions from the Heisenberg fields on ${\mathbb{R}}^3$ with co\"ordinates 
$(x,y,t)$, namely 
$$X=\partial/\partial x\qquad Y=\partial/\partial t+x\partial/\partial y
\qquad Z\equiv[X,Y]=\partial/\partial y.$$
The {\em Rumin complex\/}, discussed in Appendix~\ref{rumin}, says that 
\begin{equation}\label{Rumincomplex}
\raisebox{-30pt}{\begin{picture}(276,65)(0,-25)
\put(0,20){\makebox(0,0){$0$}}
\put(10,20){\vector(1,0){20}}
\put(40,20){\makebox(0,0){${\mathcal{E}}_3$}}
\put(50,22){\vector(2,1){20}}
\put(50,18){\vector(2,-1){20}}
\put(80,34){\makebox(0,0){${\mathcal{E}}_3$}}
\put(80,20){\makebox(0,0){$\oplus$}}
\put(80,6){\makebox(0,0){${\mathcal{E}}_3$}}
\put(90,32){\vector(4,-1){96}}
\put(90,8){\vector(4,1){96}}
\put(90,34){\vector(1,0){96}}
\put(90,6){\vector(1,0){96}}
\put(196,34){\makebox(0,0){${\mathcal{E}}_3$}}
\put(196,20){\makebox(0,0){$\oplus$}}
\put(196,6){\makebox(0,0){${\mathcal{E}}_3$}}
\put(206,32){\vector(2,-1){20}}
\put(206,8){\vector(2,1){20}}
\put(236,20){\makebox(0,0){${\mathcal{E}}_3$}}
\put(246,20){\vector(1,0){20}}
\put(276,20){\makebox(0,0){$0$}}
\put(50,-15){\makebox(0,0){\scriptsize$f\mapsto\left[\!\!\begin{array}c Xf\\
Yf\end{array}\!\!\right]$}}
\put(138,-15){\makebox(0,0){\scriptsize$\left[\!\!\begin{array}c a\\
b\end{array}\!\!\right]\mapsto\left[\!\!\begin{array}c X^2b-(XY+Z)a\\
Y^2a-(YX-Z)b\end{array}\!\!\right]$}}
\put(237,-15){\makebox(0,0){\scriptsize$\left[\!\!\begin{array}c c\\
d\end{array}\!\!\right]\mapsto Xd+Yc$}}
\end{picture}}\end{equation}
is locally exact except at the leftmost ${\mathcal{E}}_3$ where the cohomology 
is~${\mathbb{R}}$, the real-valued locally constant functions. Here, we are 
writing ${\mathcal{E}}_3$ for the germs of smooth functions of the three 
variables $(x,y,t)$ and shortly we shall write ${\mathcal{E}}_2$ for the germs 
of smooth functions of the two variables~$(x,y)$. Evidently, there is a short 
exact sequence
$$0\to{\mathcal{E}}_2\to{\mathcal{E}}_3\xrightarrow{\,\partial/\partial t\,}
{\mathcal{E}}_3\to 0.$$
Also note that the vector field $\partial/\partial t$ commutes with $X,Y,Z$. 
Therefore, we may consider the commutative diagram
$$\begin{picture}(276,110)(0,-70)
\put(0,20){\makebox(0,0){$0$}}
\put(10,20){\vector(1,0){20}}
\put(40,20){\makebox(0,0){${\mathcal{E}}_3$}}
\put(50,22){\vector(2,1){20}}
\put(50,18){\vector(2,-1){20}}
\put(80,34){\makebox(0,0){${\mathcal{E}}_3$}}
\put(80,20){\makebox(0,0){$\oplus$}}
\put(80,6){\makebox(0,0){${\mathcal{E}}_3$}}
\put(90,32){\vector(4,-1){96}}
\put(90,8){\vector(4,1){96}}
\put(90,34){\vector(1,0){96}}
\put(90,6){\vector(1,0){96}}
\put(196,34){\makebox(0,0){${\mathcal{E}}_3$}}
\put(196,20){\makebox(0,0){$\oplus$}}
\put(196,6){\makebox(0,0){${\mathcal{E}}_3$}}
\put(206,32){\vector(2,-1){20}}
\put(206,8){\vector(2,1){20}}
\put(236,20){\makebox(0,0){${\mathcal{E}}_3$}}
\put(246,20){\vector(1,0){20}}
\put(276,20){\makebox(0,0){$0$}}
\put(40,-40){\vector(0,1){50}}
\put(80,-25){\vector(0,1){20}}
\put(196,-25){\vector(0,1){20}}
\put(236,-40){\vector(0,1){50}}
\put(52,-15){\makebox(0,0){\scriptsize$\partial/\partial t$}}
\put(92,-15){\makebox(0,0){\scriptsize$\partial/\partial t$}}
\put(208,-15){\makebox(0,0){\scriptsize$\partial/\partial t$}}
\put(248,-15){\makebox(0,0){\scriptsize$\partial/\partial t$}}
\put(0,-50){\makebox(0,0){$0$}}
\put(10,-50){\vector(1,0){20}}
\put(40,-50){\makebox(0,0){${\mathcal{E}}_3$}}
\put(50,-48){\vector(2,1){20}}
\put(50,-52){\vector(2,-1){20}}
\put(80,-36){\makebox(0,0){${\mathcal{E}}_3$}}
\put(80,-50){\makebox(0,0){$\oplus$}}
\put(80,-64){\makebox(0,0){${\mathcal{E}}_3$}}
\put(90,-38){\vector(4,-1){96}}
\put(90,-62){\vector(4,1){96}}
\put(90,-36){\vector(1,0){96}}
\put(90,-64){\vector(1,0){96}}
\put(196,-36){\makebox(0,0){${\mathcal{E}}_3$}}
\put(196,-50){\makebox(0,0){$\oplus$}}
\put(196,-64){\makebox(0,0){${\mathcal{E}}_3$}}
\put(206,-38){\vector(2,-1){20}}
\put(206,-62){\vector(2,1){20}}
\put(236,-50){\makebox(0,0){${\mathcal{E}}_3$}}
\put(246,-50){\vector(1,0){20}}
\put(276,-50){\makebox(0,0){$0$}}
\end{picture}$$
to which we may apply the spectral sequences of a double complex, or indulge in
diagram chasing, to conclude that, not only is there a complex of
differential operators
$$\begin{picture}(276,65)(0,-25)
\put(0,20){\makebox(0,0){$0$}}
\put(10,20){\vector(1,0){20}}
\put(40,20){\makebox(0,0){${\mathcal{E}}_2$}}
\put(50,22){\vector(2,1){20}}
\put(50,18){\vector(2,-1){20}}
\put(80,34){\makebox(0,0){${\mathcal{E}}_2$}}
\put(80,20){\makebox(0,0){$\oplus$}}
\put(80,6){\makebox(0,0){${\mathcal{E}}_2$}}
\put(90,32){\vector(4,-1){96}}
\put(90,8){\vector(4,1){96}}
\put(90,34){\vector(1,0){96}}
\put(90,6){\vector(1,0){96}}
\put(196,34){\makebox(0,0){${\mathcal{E}}_2$}}
\put(196,20){\makebox(0,0){$\oplus$}}
\put(196,6){\makebox(0,0){${\mathcal{E}}_2$}}
\put(206,32){\vector(2,-1){20}}
\put(206,8){\vector(2,1){20}}
\put(236,20){\makebox(0,0){${\mathcal{E}}_2$}}
\put(246,20){\vector(1,0){20}}
\put(276,20){\makebox(0,0){$0$}}
\put(50,-15){\makebox(0,0){\scriptsize$f\mapsto\left[\!\!\begin{array}c Xf\\
Yf\end{array}\!\!\right]$}}
\put(138,-15){\makebox(0,0){\scriptsize$\left[\!\!\begin{array}c a\\
b\end{array}\!\!\right]\mapsto\left[\!\!\begin{array}c X^2b-(XY+Z)a\\
Y^2a-(YX-Z)b\end{array}\!\!\right]$}}
\put(237,-15){\makebox(0,0){\scriptsize$\left[\!\!\begin{array}c c\\
d\end{array}\!\!\right]\mapsto Xd+Yc$}}
\end{picture}$$
in which $X,Y,Z$ now denote the differential operators (\ref{firstgrushin})
on~${\mathbb{R}}^2$, but also that the cohomology of this complex resides in 
the zeroth and first degrees, where it is~${\mathbb{R}}$. In particular, we 
have found integrability conditions for two smooth functions $a=a(x,y)$ and 
$b=b(x,y)$ to be locally of the form $a=Xf$ and $b=Yf$ for some $f=f(x,y)$ as 
follows.
\begin{thm}
Suppose\/ $U^{\mathrm{open}}\subset{\mathbb{R}}^2$ is contractible. Then, for a
pair of smooth functions\/ $a$ and\/ $b$ defined on\/~$U$,
$$\left.\begin{array}{rcl}X^2b&=&(XY+Z)a\\
Y^2a&=&(YX-Z)b\end{array}\right\}\iff\begin{tabular}{l}
$\exists$ a smooth function\/ $f$ on\/~$U$\\
and a constant\/ $C$ such that\\
$Xf=a$ and\/ $Yf=C+b$.
\end{tabular}$$
\end{thm}

\section{The Martinet distribution}\label{Martinet}
We shall derive integrability conditions for the Martinet fields 
(\ref{martinet}) on ${\mathbb{R}}^3$ from the Engel complex on 
${\mathbb{R}}^4$ constructed from the fields
\begin{equation}\label{Engelfields}
\begin{array}{ll}X=\partial/\partial x&\quad 
Y=\partial/\partial z + x\partial/\partial t+ x^2\partial/\partial y\\
Z\equiv[X,Y]=\partial/\partial t+2x\partial/\partial y&\quad
W\equiv[X,Z]=2\partial/\partial y.
\end{array}\end{equation}
The Engel complex takes the form
$$0\to{\mathcal{E}}_4
\begin{array}{c}\nearrow\\ \searrow\end{array}\!\!
\begin{array}c{\mathcal{E}}_4\\[2pt] \oplus\\[4pt]
{\mathcal{E}}_4\end{array}\!\!
\begin{array}c\longrightarrow\\ \mbox{\Large\begin{picture}(0,0)
\put(0,0){\makebox(0,0){$\nearrow$}}
\put(0,0){\makebox(0,0){$\searrow$}}
\end{picture}}\\[4pt] \longrightarrow\end{array}\!\!
\begin{array}c{\mathcal{E}}_4\\[2pt] \oplus\\[4pt]
{\mathcal{E}}_4\end{array}\!\!
\begin{array}c\longrightarrow\\ \mbox{\Large\begin{picture}(0,0)
\put(0,0){\makebox(0,0){$\nearrow$}}
\put(0,0){\makebox(0,0){$\searrow$}}
\end{picture}}\\[4pt] \longrightarrow\end{array}\!\!
\begin{array}c{\mathcal{E}}_4\\[2pt] \oplus\\[4pt]
{\mathcal{E}}_4\end{array}\!\!
\begin{array}{c}\searrow\\ \nearrow\end{array}
{\mathcal{E}}_4\to 0,$$
where the differential operators are
\begin{equation}\label{Engeloperators}\raisebox{10pt}{\makebox[200pt]{\small
$\begin{array}{l}f\mapsto
\left[\!\!\begin{array}{c}Xf\\ Yf\end{array}\!\!\right]
\enskip
\left[\!\!\begin{array}{c}a\\ b\end{array}\!\!\right]\mapsto
\left[\!\!\begin{array}{c}X^3b-(X^2Y+XZ+W)a\\ 
Y^2a-(YX-Z)b\end{array}\!\!\right]\\[15pt]
\hspace*{80pt}
\left[\!\!\begin{array}{c}c\\ d\end{array}\!\!\right]\mapsto
\left[\!\!\begin{array}{c}X^3d+(XY+Z)c\\ 
                          Y^2c+(YX^2-ZX+W)d\end{array}\!\!\right]
\enskip\left[\!\!\begin{array}{c}g\\ h\end{array}\!\!\right]\mapsto
Xh-Yg.
\end{array}$}}\end{equation}
Notice that $\partial/\partial t$ commutes with each of the vector fields 
$X,Y,Z,W$ and so we have a commutative diagram
$$\begin{array}{ccccccccccccl}
0&\to&{\mathcal{E}}_4&\to&{\mathcal{E}}_4^2&\to&{\mathcal{E}}_4^2
&\to&{\mathcal{E}}_4^2&\to&{\mathcal{E}}_4&\to&0\\
&&\Big\uparrow\makebox[0pt][l]{\scriptsize$\partial/\partial t$}
&&\Big\uparrow\makebox[0pt][l]{\scriptsize$\partial/\partial t$}
&&\Big\uparrow\makebox[0pt][l]{\scriptsize$\partial/\partial t$}
&&\Big\uparrow\makebox[0pt][l]{\scriptsize$\partial/\partial t$}
&&\Big\uparrow\makebox[0pt][l]{\scriptsize$\partial/\partial t$}\\
0&\to&{\mathcal{E}}_4&\to&{\mathcal{E}}_4^2&\to&{\mathcal{E}}_4^2
&\to&{\mathcal{E}}_4^2&\to&{\mathcal{E}}_4&\to&0,
\end{array}$$
where each row is the Engel complex. Arguing as in \S\ref{firstGrushin}, we 
conclude that there is a complex of differential operators on ${\mathbb{R}}^3$
$$0\to{\mathcal{E}}_3
\begin{array}{c}\nearrow\\ \searrow\end{array}\!\!
\begin{array}c{\mathcal{E}}_3\\[2pt] \oplus\\[4pt]
{\mathcal{E}}_3\end{array}\!\!
\begin{array}c\longrightarrow\\ \mbox{\Large\begin{picture}(0,0)
\put(0,0){\makebox(0,0){$\nearrow$}}
\put(0,0){\makebox(0,0){$\searrow$}}
\end{picture}}\\[4pt] \longrightarrow\end{array}\!\!
\begin{array}c{\mathcal{E}}_3\\[2pt] \oplus\\[4pt]
{\mathcal{E}}_3\end{array}\!\!
\begin{array}c\longrightarrow\\ \mbox{\Large\begin{picture}(0,0)
\put(0,0){\makebox(0,0){$\nearrow$}}
\put(0,0){\makebox(0,0){$\searrow$}}
\end{picture}}\\[4pt] \longrightarrow\end{array}\!\!
\begin{array}c{\mathcal{E}}_3\\[2pt] \oplus\\[4pt]
{\mathcal{E}}_3\end{array}\!\!
\begin{array}{c}\searrow\\ \nearrow\end{array}
{\mathcal{E}}_3\to 0,$$
in which the differential operators are exactly as in (\ref{Engeloperators}) 
except that $X,Y,Z,W$ now stand for the vector fields 
$$\begin{array}{ll}X=\partial/\partial x&\quad 
Y=\partial/\partial z + x^2\partial/\partial y\\
Z\equiv[X,Y]=2x\partial/\partial y&\quad
W\equiv[X,Z]=2\partial/\partial y
\end{array}$$
on ${\mathbb{R}}^3$, the first two of which are the Martinet 
fields~(\ref{martinet}). Moreover, the local cohomology of the complex occurs 
only in degrees zero and one, where it is ${\mathbb{R}}$. In particular, the 
integrability conditions for the Martinet fields as as follows.
\begin{thm}
Suppose\/ $U^{\mathrm{open}}\subset{\mathbb{R}}^3$ is contractible. Let\/ $X$
and\/ $Y$ denote the Martinet fields\/ {\rm(\ref{martinet})}
on\/~${\mathbb{R}}^3$. Then, for a pair of smooth functions\/ $a$ and\/ $b$
defined on\/~$U$, 
$$\left.\begin{array}{rcl}X^3b&=&(X^2Y+XZ+W)a\\
Y^2a&=&(YX-Z)b\end{array}\right\}\!\iff\!\begin{tabular}{l}
$\exists$ a smooth function\/ $f$ on\/~$U$\\
and a constant\/ $C$ such that\\
$Xf=a$ and\/ $Yf=Cx+b$.
\end{tabular}$$
Here, for convenience, we have set\/ $Z\equiv[X,Y]$ and\/~$W\equiv[X,Z]$.
\end{thm}

\section{The second Grushin distribution}
Recall that we are concerned with the four vector fields
\begin{equation}\label{secondGrushinfields}
X=\partial/\partial x\qquad Y=x^2\partial/\partial y\qquad
Z=2x\partial/\partial y\qquad W=2\partial/\partial y\end{equation}
on ${\mathbb{R}}^2$, which we may clearly view as the four Engel fields 
(\ref{Engelfields}) acting on smooth functions $f$ on ${\mathbb{R}}^4$ that 
happen to be of the form~$f=f(x,y)$. Evidently, we have the exact sequence
\begin{equation}\label{rel_deRham}\begin{picture}(200,40)
\put(0,20){\makebox(0,0){$\;0$}}
\put(10,20){\vector(1,0){20}}
\put(40,20){\makebox(0,0){${\mathcal{E}}_2$}}
\put(50,20){\vector(1,0){20}}
\put(80,20){\makebox(0,0){${\mathcal{E}}_4$}}
\put(90,22){\vector(2,1){20}}
\put(90,18){\vector(2,-1){20}}
\put(96,34){\makebox(0,0){\scriptsize$\partial/\partial z$}}
\put(96,6){\makebox(0,0){\scriptsize$\partial/\partial t$}}
\put(144,34){\makebox(0,0){\scriptsize$\partial/\partial t$}}
\put(144,6){\makebox(0,0){\scriptsize$-\partial/\partial z$}}
\put(120,34){\makebox(0,0){${\mathcal{E}}_4$}}
\put(120,20){\makebox(0,0){$\oplus$}}
\put(120,6){\makebox(0,0){${\mathcal{E}}_4$}}
\put(130,32){\vector(2,-1){20}}
\put(130,8){\vector(2,1){20}}
\put(160,20){\makebox(0,0){${\mathcal{E}}_4$}}
\put(170,20){\vector(1,0){20}}
\put(200,20){\makebox(0,0){$0$.}}
\end{picture}\end{equation}
It is the de~Rham complex in the $(z,t)$-variables. Since both 
$\partial/\partial z$ and $\partial/\partial t$ commute with the Engel fields 
(\ref{Engelfields}) there is a commutative diagram 
$$\addtolength{\arraycolsep}{-2pt}\begin{array}{ccccccccccccl}
0&\to&{\mathcal{E}}_4&\to&{\mathcal{E}}_4^2&\to&{\mathcal{E}}_4^2
&\to&{\mathcal{E}}_4^2&\to&{\mathcal{E}}_4&\to&0\\
&&\nearrow\;\nwarrow
&&\nearrow\;\nwarrow
&&\nearrow\;\nwarrow
&&\nearrow\;\nwarrow
&&\nearrow\;\nwarrow\\
0&\to&{\mathcal{E}}_4\oplus{\mathcal{E}}_4&\to&
{\mathcal{E}}_4^2\oplus{\mathcal{E}}_4^2&\to&
{\mathcal{E}}_4^2\oplus{\mathcal{E}}_4^2&\to&
{\mathcal{E}}_4^2\oplus{\mathcal{E}}_4^2&\to&
{\mathcal{E}}_4\oplus{\mathcal{E}}_4&\to&0\\
&&\nwarrow\;\nearrow
&&\nwarrow\;\nearrow
&&\nwarrow\;\nearrow
&&\nwarrow\;\nearrow
&&\nwarrow\;\nearrow\\
0&\to&{\mathcal{E}}_4&\to&{\mathcal{E}}_4^2&\to&{\mathcal{E}}_4^2
&\to&{\mathcal{E}}_4^2&\to&{\mathcal{E}}_4&\to&0,
\end{array}$$
in which the rows are the Engel complex with appropriate multiplicity and the
columns are copies of (\ref{rel_deRham}). As in \S\ref{firstGrushin}
and~\S\ref{Martinet}, it follows that there is a complex of differential 
operators on~${\mathbb{R}}^2$
\begin{equation}\label{complexforGrushin2}
0\to{\mathcal{E}}_2
\begin{array}{c}\nearrow\\ \searrow\end{array}\!\!
\begin{array}c{\mathcal{E}}_2\\[2pt] \oplus\\[4pt]
{\mathcal{E}}_2\end{array}\!\!
\begin{array}c\longrightarrow\\ \mbox{\Large\begin{picture}(0,0)
\put(0,0){\makebox(0,0){$\nearrow$}}
\put(0,0){\makebox(0,0){$\searrow$}}
\end{picture}}\\[4pt] \longrightarrow\end{array}\!\!
\begin{array}c{\mathcal{E}}_2\\[2pt] \oplus\\[4pt]
{\mathcal{E}}_2\end{array}\!\!
\begin{array}c\longrightarrow\\ \mbox{\Large\begin{picture}(0,0)
\put(0,0){\makebox(0,0){$\nearrow$}}
\put(0,0){\makebox(0,0){$\searrow$}}
\end{picture}}\\[4pt] \longrightarrow\end{array}\!\!
\begin{array}c{\mathcal{E}}_2\\[2pt] \oplus\\[4pt]
{\mathcal{E}}_2\end{array}\!\!
\begin{array}{c}\searrow\\ \nearrow\end{array}
{\mathcal{E}}_2\to 0,\end{equation}
in which the differential operators are exactly as in (\ref{Engeloperators}) 
except that $X,Y,Z,W$ now stand for the vector 
fields~(\ref{secondGrushinfields}). Moreover, this complex has local 
cohomology ${\mathbb{R}}$ in degree zero, ${\mathbb{R}}\oplus{\mathbb{R}}$ in 
degree one, and ${\mathbb{R}}$ in degree two. More explicitly, we have proved 
the following.
\begin{thm}
Suppose\/ $U^{\mathrm{open}}\subset{\mathbb{R}}^2$ is contractible and let\/ 
$X,Y,Z,W$ denote the Grushin fields\/ {\rm(\ref{secondGrushinfields})}
on\/~${\mathbb{R}}^2$. Then, for smooth functions\/ $a$ and\/ $b$
defined on\/~$U$, 
$$\begin{array}{rcl}X^3b&=&(X^2Y+XZ+W)a\\
Y^2a&=&(YX-Z)b\end{array}$$
if and only if there is a smooth function\/ $f$ on\/ $U$ and constants\/ $C$
and\/~$D$ such that
$Xf=a$ and
$Yf=Cx+D+b$.
Furthermore, for smooth functions\/ $c$ and\/ $d$ on\/~$U$,
$$\begin{array}{rcl}X^3d+(XY+Z)c&=&0\\
Y^2c+(YX^2-ZX+W)d&=&0\end{array}$$
if and only if there are smooth functions\/ $a$ and\/ $b$ on\/ $U$
and a constant\/ $E$ such that
$$\begin{array}{rcl}X^3b-(X^2Y+XZ+W)a&=&c\\
Y^2a-(YX-Z)b&=&E+d.
\end{array}$$
Otherwise, the complex\/ {\rm(\ref{complexforGrushin2})} is locally exact.
\end{thm}

\appendix
\section{The Rumin complex on~${\mathbb{R}}^3$}\label{rumin}
There are several different viewpoints on the Rumin complex on ${\mathbb{R}}^3$
and here is not the place to go into the details concerning various subtle
distinctions. In fact, the minimal structure required for the basic
construction is that of a contact distribution~\cite{R}. A more refined outcome
is obtained starting with a pair of line fields that together span a contact
distribution. This structure is known as {\em contact-Lagrangian\/}
\cite[\S4.2.3]{CSbook} or sometimes as {\em para-CR\/}~\cite{AMT}. For our
purposes, it will suffice to consider the so-called {\em flat model\/}, which
may be defined as follows. Choose vector fields $X$ and $Y$ spanning the two
line fields and let $Z\equiv[X,Y]$. The contact condition is precisely that
$X,Y,Z$ be linearly independent. For the flat model we require that $X$ and $Y$
can be chosen so that $[X,Z]$ and $[Y,Z]$ both vanish (in general, there is a
curvature obstruction to this being possible).
Following~\cite[\S8.1]{Beastwood}, the Rumin complex in this case can be
constructed from the de~Rham sequence as follows. If we denote by 
$\xi,\eta,\zeta$ the co-frame dual to $X,Y,Z$, then 
\begin{equation}\label{EDS}
d\zeta=\eta\wedge\xi\qquad d\xi=0\qquad d\eta=0\end{equation}
and we may contemplate the de~Rham complex written with respect to this 
co-frame. Specifically, let us consider the diagram
$$\begin{array}{ccccccccccc}
&&&&0&&0\\
&&&&\uparrow&&\uparrow\\
&&&&\langle\xi,\eta\rangle&&\langle\eta\wedge\xi\rangle
&&\makebox[20pt]{$\langle\eta\wedge\xi\wedge\zeta\rangle$}\\
&&&&\uparrow&&\uparrow&&\|\\
0&\to&\Lambda^0&\xrightarrow{\,d\,}&\Lambda^1&\xrightarrow{\,d\,}
&\Lambda^2&\xrightarrow{\,d\,}&\Lambda^3&\to&0,\\
&&&&\uparrow&&\uparrow\\
&&&&\langle\zeta\rangle
&&\makebox[0pt]{$\langle\xi\wedge\zeta,\zeta\wedge\eta\rangle$}\\
&&&&\uparrow&&\uparrow\\
&&&&0&&0
\end{array}$$
where $\langle\underbar\quad\rangle$ denotes the bundle spanned by the 
enclosed forms. Notice that the composition
$$\langle\zeta\rangle\to\Lambda^1\xrightarrow{\,d\,}\Lambda^2\to
\langle\eta\wedge\xi\rangle$$
is simply
$$g\,\zeta\mapsto d(g\,\zeta)=dg\wedge\zeta+g\,d\zeta=
dg\wedge\zeta+g\,\eta\wedge\xi\mapsto g\,\eta\wedge\xi$$
and hence defines an isomorphism between these line bundles. The Rumin complex
is obtained by using this isomorphism to cancel these line bundles hence
obtaining, by dint of diagram chasing, a new locally exact complex
\begin{equation}\label{prototypeBGG}0\to\Lambda^0
\begin{array}{c}\nearrow\\ \searrow\end{array}\!\!
\begin{array}c\langle\xi\rangle\\[2pt] \oplus\\[4pt]
\langle\eta\rangle\end{array}\!\!
\begin{array}c\longrightarrow\\ \mbox{\Large\begin{picture}(0,0)
\put(0,0){\makebox(0,0){$\nearrow$}}
\put(0,0){\makebox(0,0){$\searrow$}}
\end{picture}}\\[4pt] \longrightarrow\end{array}\!\!
\begin{array}c\langle\xi\wedge\zeta\rangle\\[2pt] \oplus\\[4pt]
\langle\zeta\wedge\eta\rangle\end{array}\!\!
\begin{array}{c}\searrow\\ \nearrow\end{array}
\langle\eta\wedge\xi\wedge\zeta\rangle\to 0.\end{equation}
The operators in this complex may be explicitly computed. The one-form 
$\omega=a\xi+b\eta$ for example, enjoys a unique lift $\tilde\omega$ 
annihilated by the composition 
$\Lambda^1\xrightarrow{\,d\,}\Lambda^2\to\langle\eta\wedge\xi\rangle$. 
Specifically, from (\ref{EDS}) we see that
$$\begin{array}{l}d(a\xi+b\eta+(Xb-Ya)\zeta)=\\[4pt]
\qquad(X(Xb-Ya)-Za)\xi\wedge\zeta
+(Y(Ya-Xb)+Zb)\zeta\wedge\eta
\end{array}$$
and the formul{\ae} of (\ref{Rumincomplex}) emerge. 

In fact, as explained in~\cite{Beastwood,CSbook}, the flat model may be
identified with the homogeneous space ${\mathrm{SL}}(3,{\mathbb{R}})/B$ where
$B$ is the subgroup consisting of upper triangular matrices and then the Rumin
complex (\ref{prototypeBGG}) is more accurately identified as a 
{\em Bernstein-Gelfand-Gelfand (BGG)\/} complex
$$0\to\begin{picture}(24,5)
\put(4,1.5){\line(1,0){16}}
\put(4,1.5){\makebox(0,0){$\times$}}
\put(20,1.5){\makebox(0,0){$\times$}}
\put(4,8){\makebox(0,0){$\scriptstyle 0$}}
\put(20,8){\makebox(0,0){$\scriptstyle 0$}}
\end{picture}\begin{array}{c}\nearrow\\ \searrow\end{array}\!\!
\begin{array}c\begin{picture}(24,5)
\put(4,1.5){\line(1,0){16}}
\put(4,1.5){\makebox(0,0){$\times$}}
\put(20,1.5){\makebox(0,0){$\times$}}
\put(4,7.5){\makebox(0,0){$\scriptstyle -2$}}
\put(20,8){\makebox(0,0){$\scriptstyle 1$}}
\end{picture}\\[2pt] \oplus\\[4pt]
\begin{picture}(24,5)
\put(4,1.5){\line(1,0){16}}
\put(4,1.5){\makebox(0,0){$\times$}}
\put(20,1.5){\makebox(0,0){$\times$}}
\put(4,8){\makebox(0,0){$\scriptstyle 1$}}
\put(20,7.5){\makebox(0,0){$\scriptstyle -2$}}
\end{picture}\end{array}\!\!
\begin{array}c\longrightarrow\\ \mbox{\Large\begin{picture}(0,0)
\put(0,0){\makebox(0,0){$\nearrow$}}
\put(0,0){\makebox(0,0){$\searrow$}}
\end{picture}}\\[4pt] \longrightarrow\end{array}\!\!
\begin{array}c\begin{picture}(24,5)
\put(4,1.5){\line(1,0){16}}
\put(4,1.5){\makebox(0,0){$\times$}}
\put(20,1.5){\makebox(0,0){$\times$}}
\put(4,7.5){\makebox(0,0){$\scriptstyle -3$}}
\put(20,8){\makebox(0,0){$\scriptstyle 0$}}
\end{picture}\\[2pt] \oplus\\[4pt]
\begin{picture}(24,5)
\put(4,1.5){\line(1,0){16}}
\put(4,1.5){\makebox(0,0){$\times$}}
\put(20,1.5){\makebox(0,0){$\times$}}
\put(4,8){\makebox(0,0){$\scriptstyle 0$}}
\put(20,7.5){\makebox(0,0){$\scriptstyle -3$}}
\end{picture}\end{array}\!\!
\begin{array}{c}\searrow\\ \nearrow\end{array}
\begin{picture}(24,5)
\put(4,1.5){\line(1,0){16}}
\put(4,1.5){\makebox(0,0){$\times$}}
\put(20,1.5){\makebox(0,0){$\times$}}
\put(4,7.5){\makebox(0,0){$\scriptstyle -2$}}
\put(20,7.5){\makebox(0,0){$\scriptstyle -2$}}
\end{picture}\to 0.$$

\section{The Engel complex on~${\mathbb{R}}^4$}\label{engel}
This case is discussed in detail in~\cite[\S\S3,7]{BEGN}. Here, suffice it to
say that the construction is completely parallel to that just discussed and
that the result is a BGG complex for the homogeneous space
${\mathrm{Sp}}(4,{\mathbb{R}})/B$ with $B$ a Borel subgroup. In the notation
of~\cite{Beastwood} we obtain
$$0\to\begin{picture}(24,5)
\put(5.6,0){\line(1,0){12.8}}
\put(5.6,3.2){\line(1,0){12.8}}
\put(4,1.5){\makebox(0,0){$\times$}}
\put(20,1.5){\makebox(0,0){$\times$}}
\put(12,1.5){\makebox(0,0){$\langle$}}
\put(4,8){\makebox(0,0){$\scriptstyle 0$}}
\put(20,8){\makebox(0,0){$\scriptstyle 0$}}
\end{picture}\begin{array}{c}\nearrow\\ \searrow\end{array}\!\!
\begin{array}c\begin{picture}(24,5)
\put(5.6,0){\line(1,0){12.8}}
\put(5.6,3.2){\line(1,0){12.8}}
\put(4,1.5){\makebox(0,0){$\times$}}
\put(20,1.5){\makebox(0,0){$\times$}}
\put(12,1.5){\makebox(0,0){$\langle$}}
\put(4,7.5){\makebox(0,0){$\scriptstyle -2$}}
\put(20,8){\makebox(0,0){$\scriptstyle 1$}}
\end{picture}\\[2pt] \oplus\\[4pt]
\begin{picture}(24,5)
\put(5.6,0){\line(1,0){12.8}}
\put(5.6,3.2){\line(1,0){12.8}}
\put(4,1.5){\makebox(0,0){$\times$}}
\put(20,1.5){\makebox(0,0){$\times$}}
\put(12,1.5){\makebox(0,0){$\langle$}}
\put(4,8){\makebox(0,0){$\scriptstyle 2$}}
\put(20,7.5){\makebox(0,0){$\scriptstyle -2$}}
\end{picture}\end{array}\!\!
\begin{array}c\longrightarrow\\ \mbox{\Large\begin{picture}(0,0)
\put(0,0){\makebox(0,0){$\nearrow$}}
\put(0,0){\makebox(0,0){$\searrow$}}
\end{picture}}\\[4pt] \longrightarrow\end{array}\!\!
\begin{array}c\begin{picture}(24,5)
\put(5.6,0){\line(1,0){12.8}}
\put(5.6,3.2){\line(1,0){12.8}}
\put(4,1.5){\makebox(0,0){$\times$}}
\put(20,1.5){\makebox(0,0){$\times$}}
\put(12,1.5){\makebox(0,0){$\langle$}}
\put(4,7.5){\makebox(0,0){$\scriptstyle -4$}}
\put(20,8){\makebox(0,0){$\scriptstyle 1$}}
\end{picture}\\[2pt] \oplus\\[4pt]
\begin{picture}(24,5)
\put(5.6,0){\line(1,0){12.8}}
\put(5.6,3.2){\line(1,0){12.8}}
\put(4,1.5){\makebox(0,0){$\times$}}
\put(20,1.5){\makebox(0,0){$\times$}}
\put(12,1.5){\makebox(0,0){$\langle$}}
\put(4,8){\makebox(0,0){$\scriptstyle 2$}}
\put(20,7.5){\makebox(0,0){$\scriptstyle -3$}}
\end{picture}\end{array}\!\!
\begin{array}c\longrightarrow\\ \mbox{\Large\begin{picture}(0,0)
\put(0,0){\makebox(0,0){$\nearrow$}}
\put(0,0){\makebox(0,0){$\searrow$}}
\end{picture}}\\[4pt] \longrightarrow\end{array}\!\!
\begin{array}c\begin{picture}(24,5)
\put(5.6,0){\line(1,0){12.8}}
\put(5.6,3.2){\line(1,0){12.8}}
\put(4,1.5){\makebox(0,0){$\times$}}
\put(20,1.5){\makebox(0,0){$\times$}}
\put(12,1.5){\makebox(0,0){$\langle$}}
\put(4,7.5){\makebox(0,0){$\scriptstyle -4$}}
\put(20,8){\makebox(0,0){$\scriptstyle 0$}}
\end{picture}\\[2pt] \oplus\\[4pt]
\begin{picture}(24,5)
\put(5.6,0){\line(1,0){12.8}}
\put(5.6,3.2){\line(1,0){12.8}}
\put(4,1.5){\makebox(0,0){$\times$}}
\put(20,1.5){\makebox(0,0){$\times$}}
\put(12,1.5){\makebox(0,0){$\langle$}}
\put(4,8){\makebox(0,0){$\scriptstyle 0$}}
\put(20,7.5){\makebox(0,0){$\scriptstyle -3$}}
\end{picture}\end{array}\!\!
\begin{array}{c}\searrow\\ \nearrow\end{array}
\begin{picture}(24,5)
\put(5.6,0){\line(1,0){12.8}}
\put(5.6,3.2){\line(1,0){12.8}}
\put(4,1.5){\makebox(0,0){$\times$}}
\put(20,1.5){\makebox(0,0){$\times$}}
\put(12,1.5){\makebox(0,0){$\langle$}}
\put(4,7.5){\makebox(0,0){$\scriptstyle -2$}}
\put(20,7.5){\makebox(0,0){$\scriptstyle -2$}}
\end{picture}\to 0$$
and following the construction in \cite{BEGN} with a co-frame
$\xi,\eta,\zeta,\omega$ such that 
$$d\omega=\zeta\wedge\xi\qquad
d\zeta=\eta\wedge\xi\qquad d\xi=0\qquad
d\eta=0$$
leads to the explicit formul{\ae} of~(\ref{Engeloperators}). We remark that 
the Engel fields are often written as
$$X=\partial/\partial x-z\partial/\partial t-t\partial/\partial y\quad 
Y=\partial/\partial z
\quad Z=\partial/\partial t\quad
W=\partial/\partial y$$
but we prefer the form (\ref{Engelfields}) so as better to relate to the
Martinet and Grushin distributions.


\begin{thebibliography}{11}
    
\bibitem{AMT} D.V. Alekseevsky, C. Medori, and A. Tomassini,
{\em Maximally homogeneous para-CR manifolds},
Ann. Global Anal. Geom. {\bf 30} (2006), 1--27.

\bibitem{Beastwood} R.J. Baston and M.G. Eastwood,
{\em The Penrose Transform: its Interaction with Representation Theory},
Oxford University Press 1989.

\bibitem{BEGN} R.L. Bryant, M.G. Eastwood, A.R. Gover, and K. Neusser
{\em Some differential complexes within and beyond parabolic geometry},
arXiv:1112.2142.

\bibitem{CSbook} A. \v{C}ap and J. Slov\'ak,
{\em Parabolic Geometries I: Background and General Theory}, 
Amer. Math. Soc. 2009.

\bibitem{G} V.V. Grushin,
{\em A certain class of hypoelliptic operators},
Mat. Sb. {\bf 83} (1970) 456--473.

\bibitem{M} J. Martinet,
{\em Sur les singularit\'es des formes diff\'erentielles},
Ann. Inst. Fourier {\bf 20} (1970) 90--178.

\bibitem{R} M. Rumin,
{\em Un complexe de formes diff\'erentielles sur les vari\'et\'es de contact},
Comptes Rendus Acad. Sci. Paris Math. {\bf 310} (1990) 401--404.

\end{thebibliography}
\end{document}